\documentclass[10pt]{article}
\usepackage{amsmath,amsthm,amssymb}
\usepackage[dvipdfmx]{graphicx}
\usepackage{color}
\usepackage{wrapfig}
\newcommand{\maprightu}[1]{%
\smash{\mathop{%
\hbox to 1cm{\rightarrowfill}}\limits^{#1}}}
\newcommand{\maprightd}[1]{%
\smash{\mathop{%
\hbox to 1cm{\rightarrowfill}}\limits_{#1}}}
\newcommand{\mapleftu}[1]{%
\smash{\mathop{%
\hbox to 1cm{\leftarrowfill}}\limits^{#1}}}
\newcommand{\mapleftd}[1]{%
\smash{\mathop{%
\hbox to 1cm{\leftarrowfill}}\limits_{#1}}}

\newcommand{\mapdownr}[1]{\Big\downarrow
\rlap{$\vcenter{\hbox{$\scriptstyle#1\,$}}$ }}
\newcommand{\mapupl}[1]{\Big\uparrow
\llap{$\vcenter{\hbox{$\scriptstyle#1\,$}}$ }}

\begin{document}

\baselineskip 15pt

\begin{center}
{\Large \bf On Incompleteness of Some Integrable Rational Maps}
\vspace{0.5cm}

by

\vspace{0.2cm}
S.Saito$^*$, N.Saitoh$^{**}$, T.Hatanaka$^\dagger$, Y.Wakimoto$^{\dagger\dagger}$ and T.Yumibayashi$^{\dagger\dagger}$
\vspace{0.5cm}

\begin{minipage}{16cm}
{\scriptsize
$^*$ Hakusan 4-19-10, Midori-ku, Yokohama, 226-0006 Japan

$^{**}$ Department of Applied Mathematics, Yokohama National University,
Hodogaya-ku, Yokohama, 240-8501 Japan

$^\dagger$ Wakamiyacho 2-2-21, Sakata, Yamagata, 998-0053 Japan 

$^{\dagger\dagger}$ Department of Physics, Tokyo Metropolitan University,
Tokyo, 192-0397 Japan
}
\end{minipage}
\vspace{1cm}

\begin{minipage}{14cm}
{\bf Abstract:} An argument is given to associate integrable nonintegrable transition of discrete maps with the transition of Lawvere's fixed point theorem to its own contrapositive. We show that the classical description of nonlinear maps is neither complete nor totally predictable.

\end{minipage}
\end{center}

\section{Introduction}

In classical physics a dynamical system is said integrable if the states of the system are predictable precisely in all future for arbitrary initial states, otherwise it is said nonintegrable. It is well known that linear dynamical equations are integrable in general, while almost all nonlinear equations are nonintegrable. In the laws of nature the self-reference is inevitable because every small chain of action and reaction gives rise to nonlinearity of a phenomenon in physics.

The self-referential nature of nonlinearity is the source of chaotic behavior of nonintegrable systems. If a nonlinear equation is given arbitrary the probability that it is integrable is zero. Nevertheless there exist infinitely many integrable nonlinear equations, such as soliton equations, which play a fundamental role in physics. Therefore it is certainly important to know the way to characterize integrable equations among nonlinear equations. Only a little has been known about this question.

On the other hand, there have been known many paradoxes and/or incompleteness arguments in logic and mathematics, when we consider self-referential propositions. These problems, however, have never been  discussed in physics, although the same logic and mathematics are used to describe nonlinear phenomena. This is probably because the paradox arguments are thought of the subjects in logic and/or mathematics, but nothing to do with real physical phenomena. \\

The purpose of this paper is to fill this gap and explain the integrable nonintegrable transition of discrete dynamical maps by means of the fixed point theorem in logic and/or mathematics. Among various self-referential arguments, the fixed point theorem by Lawvere provides a universal description of them, such as Liar paradox, Russell's paradox, as well as Cantor's diagonal argument, G\"odel's first incompleteness theorem, Turing's halting problem, $\lambda$ fixed point functions and so on\cite{Lawvere, Yanofsky}. We show that the same type of fixed point appears at the transition between integrable and nonintegrable maps and determines their critical nature. From this fact we will conclude that the classical description of nonlinear maps is neither complete nor totally predictable.\\

We consider in this paper some higher dimensional rational maps and study the behavior of periodic points in detail. Each map has a parameter, say $a$, which can vary continuously. The map is integrable at $a=0$, but nonintegrable otherwise. This means that we consider a phase transition of the map from a nonintegrable phase to an integrable phase as $a$ changes to zero.

When the map is nonintegrable unstable periodic points form a fractal set. The closure of the set is called the Julia set\cite{Devaney}. As $a$ approaches zero the Julia set disappears altogether and periodic points form algebraic varieties, which we call invariant varieties of periodic points, or IVPPs, different for each period\cite{SS1, IVPP, SS2, SS3, SS4}. There are infinitely many IVPPs for each integrable map. Therefore the Julia set characterizes nonintegrable phase, while an IVPP characterizes integrable phase.

It is remarkable that IVPPs of all periods intersect at certain points. We call such an intersection a variety of singular points, or VSP, since every point in this variety is a point of many different periods simultaneously. Moreover we have found recently that all IVPPs are generated from the VSP by an iteration of the map\cite{SS5, YSW}. This means that the singular point is a `fixed point' of the map. We have shown in \cite{SJNMP, SYW} that this phenomenon can be associated with the projective resolution of a triangulated category. 

We are going to explain this sudden appearance of the fixed points by means of the fixed point theorem. 
Namely we associate the integrable nonintegrable phase transition with the transition of Lawvere's theorem to its contrapositive. Therefore, before studying rational maps in detail in the next section, let us review briefly the fixed point theorems. 
Lawvere's fixed point theorem states the following:
\bigskip

\noindent
{\bf Theorem (Lawvere)}

{\it In any cartesian closed category, if there exists an object $T$ and a weakly point-surjective morphism $T\to Y^T$, then $Y$ has the fixed point property.}\\

Instead of this original statement based on the category theory, however, it is more convenient for our purpose to refer the version rephrased by Yanofsky\cite{Yanofsky}.  
We say, for any sets $T$ and $Y$ and a function $\phi:  T\times T\to Y$,\\

\noindent
{\bf Definition}

$\bullet$ a function $\alpha: Y\to Y$ has a fixed point iff
 there exists $y\in Y$ satisfying $\alpha(y)=y$,

$\bullet$ a function $\psi: T\to Y$ is representable by $\phi$ iff there exists $t\in T$ satisfying $\psi(-)=\phi(-,t)$. \\

Then the following two assertions are equivalent to Lawvere's theorem.\\

\noindent
{\bf Theorem (Yanofsky version)}
\begin{enumerate}
\item[LT1] {\it If there exists a function $\alpha: Y\to Y$ which has no fixed point, there exists a function $\psi: T\to Y$ which is nonrepresentable by $\phi$.}

\item[LT2] {\it If all functions $\psi: T\to Y$ are representable by any $\phi$, all functions $\alpha: Y\to Y$ have a fixed point. }
\end{enumerate}

It will be convenient to draw them in a diagram following to the paper\cite{Yanofsky}.
\begin{equation}
\begin{array}{ccc}
T\times T&\maprightu{\phi}&Y\cr
\mapupl{\Delta}&&\mapdownr{\alpha}\cr
T&\maprightd{\psi}&Y\cr
\end{array}
\label{Yanofsky}
\end{equation}

The first statement LT1 is essentially the Cantor diagonal argument, while G\"odel's first incompleteness theorem follows to LT2, for example. LT2 is the contrapositive of LT1. Let us write this fact schematically as
\begin{equation}
\begin{array}{cccc}
{\rm LT1}:\quad&\nexists\ {\rm fixed}\ {\rm point}&\Rightarrow&
\exists\ {\rm nonrepresentable}\ {\rm function},\\
&\updownarrow&&
\updownarrow\\
{\rm LT2}:\quad&\exists\ {\rm fixed}\ {\rm point}&\Leftarrow&
\nexists\ {\rm nonrepresentable}\ {\rm function}.
\end{array}
\label{LT scheme}
\end{equation}
Here $A\Rightarrow B$ means `$A$ implies $B$', whereas $\updownarrow$ means a transition to the mutually negative statements. 
In order to compare this scheme with the transition of a nonintegrable map to integrable one, we must recall another important fact. As we show in the next section, the periodic points of all periods form either IVPPs or discrete sets. An IVPP and a set of discrete periodic points cannot exist in one map simultaneously. Let us call the discrete set of all periodic points as DSPP. {Since the Julia set is a closure of unstable periodic points of the map, we also define $\overline{{\rm DSPP}}$, the closure of DSPP, so that the Julia set is a subset of $\overline{{\rm DSPP}}$. Then} the scheme which we compare with (\ref{LT scheme}) is the following:
\begin{equation}
\begin{array}{cccccc}
{\rm nonintegrable\ map}:\quad&\nexists\ {\rm fixed\ point}&\Rightarrow&
\exists\ \overline{{\rm DSPP}}&\Rightarrow& \nexists\ {\rm IVPP},\\
&\updownarrow&&
\updownarrow&&
\updownarrow\\
{\rm integrable\ map}:\quad&\exists\ {\rm fixed\ point}&\Leftarrow&\ 
\nexists\ \overline{{\rm DSPP}}&\Leftarrow& \exists\ {\rm IVPP}.
\end{array}
\label{contrapositive}
\end{equation}
Namely the existence of a fixed point protects integrability against the Julia set, the source of chaos.

There is, however, a little gap if we want to apply LT2 directly to our problem. In Lawvere's theorem the cartesian closed category is assumed in general, whereas our objects, {\it i.e.,} IVPPs, are 
%topological
{algebraic varieties, hence it is not clear how to fill the conditions of the theorem}. In order to fill the gap we recall the 
fixed point theorem by Brouwer\cite{Brouwer}, which asserts the following:\\

\noindent
{\bf Theorem (Brouwer)}

{\it Every closed disk $D$ of any dimension has a fixed point with respect to {any} continuous map.}
\bigskip

We can associate the continuous maps in Brouwer's theorem with the generation of IVPPs of different periods in the projective resolution of the triangulated category, which we mentioned above. If we want to reinterpret this fact by means of the proposition LT2, however, we must show the correspondence between Brouwer's theorem and LT2. We will discuss this problem briefly in the last section, but leaving mathematical proof to our forthcoming paper.

In \S2 we study the behavior of periodic points of rational maps to explore the phenomenon which takes place at the transition point from nonintegrable phase to integrable one. Although our arguments are quite general we often use simple examples to illustrate our results visually. An interpretation of this phenomonon of transition to the  scheme of the fixed point theorems will be discussed in \S3. Finally the correspondence between two schemes (\ref{LT scheme}) and (\ref{contrapositive}) will be discussed in \S4.

%%%%%%%%%%%%%%%%%%%%%%%%%%%%%%%%%%%%%%%%%%%%%
\section{Integrable nonintegrable transition}

The nature of periodic points is totally different between integrable and nonintegrable maps. We discuss in this section how the transition takes place when a nonintegrable map becomes integrable . 

\subsection{The IVPP theorem}

We consider a rational map of $d$ dimension:
\[
f: x=(x_1,x_2,...,x_d)\quad \to\quad x^{(1)}=(x^{(1)}_1,x^{(1)}_2,...,x^{(1)}_d),\qquad 
x,x^{(1)}\in \mathbb{C}^d\cup\{\infty\},
\]
which has $p$ invariants $\{h_1(x),h_2(x),...,h_p(x)\}$ satisfying
\[
h_i(x)=h_i(x^{(1)}(x)),\qquad i=1,2,...,p.
\]
The periodic points of period $n$ are those satisfying the conditions
\begin{equation}
x_j^{(n)}(x):=f\overbrace{(f(f(\cdots (f(}^{n}x))\cdots )))_j=x_j,\qquad j=1,2,...,d,\quad n=2,3,4,....
\label{ppc}
\end{equation}
If the values of the invariants are fixed at $h=(h_1,h_2,...,h_p)$, the map is constrained on the level set specified by $h$, and only $d-p$ conditions of (\ref{ppc}) are independent, which we write as
\begin{equation}
\Gamma_\alpha^{(n)}(\xi,h)=0,\qquad \alpha=1,2,...,d-p,\quad n=2,3,4,.....
\label{Gamma}
\end{equation}
Here $\xi=(\xi_1,\xi_2,...,\xi_{d-p})$ is the set of $d-p$ independent  variables which specify the coordinates of the level set. Since the number of equations in (\ref{Gamma}) is the same with the number of the variables $\xi$, we find a set of discrete points on the level set. They are the set of periodic points of period $n$. 

There is, however, a possibility that the level set is filled by the periodic points of one particular period, say of period $k$. It can happen if the conditions (\ref{Gamma}) do not depend on $\xi$, but are determined by the invariants alone. In such a case we write (\ref{Gamma}) as
\begin{equation}
\gamma^{(k)}_\alpha(h)=0,\qquad \alpha=1,2,...,d-p,\quad k=2,3,4,....,
\label{gamma}
\end{equation}
which has a solution for $h$ only if $p\ge d/2$. We notice that, if the periodicity condition is given by (\ref{gamma}), all periodic points of period $k$ form an algebraic variety
\begin{equation}
v^{(k)}=\Big\{x\in\mathbb{C}^d\Big|\gamma^{(k)}_\alpha(h(x))=0,\ \alpha=1,2,...,d-p\Big\},\quad k=2,3,4,.....
\label{v^k}
\end{equation}
The dimension of this variety is $p$. Since it is determined by the invariants alone, we call $v^{(k)}$ the invariant variety of periodic points (IVPP) of period $k$. Our theorem follows immediately from the definition of the IVPP (\ref{v^k}):
\bigskip

\noindent
{\bf IVPP Theorem}\cite{SS1, IVPP, SS2, SS3, SS4}

{\it If $p\ge d/2$, an IVPP and a discrete set of periodic points of any period {on a level set} are not in one map at the same time.}\\

\noindent
Proof:\quad Suppose the period $k$ and period $n(\ne k)$ conditions are given by (\ref{gamma}) and (\ref{Gamma}), respectively. Since these two algebraic relations can be solved for $(\xi,h)$ simultaneously, the solutions are periodic points of both periods $k$ and $n$, which is impossible, in general, as far as $n\ne k$.\\

Some remarks are in order.
\begin{enumerate}
\item If there is an IVPP of any period the map has no Julia set.
\item If the existence of the Julia set is necessary for a rational map being nonintegrable, the existence of an IVPP of any period guarantees integrability of the map.
\item We have investigated many integrable rational maps satisfying $p\ge d/2$, and found IVPPs of all peiods for each map iteratively by a projective resolution of the singularity\cite{SS5, YSW, SJNMP, SYW}.
\end{enumerate}

\subsection{The phase transition}

In order to study explicitly the transition between integrable and nonintegrable phases we consider a rational map $f_a$ with a control parameter $a$. The map is nonintegrable when $a\ne 0$, while it is integrable at $a=0$.  

We want to know the behavior of periodic points of the map as $a$ changes. They form a set of discrete points DSPP when $a\ne 0$. In particular unstable periodic points together with their closure form a Julia set $J_a$. If $p_0$ is any unstable periodic point the Julia set can be obtained by repeating the inverse map $f_a^{-1}$ iteratively starting from $p_0$. 
\begin{equation}
J_a=\overline{\left\{\bigcup_{k=0}^\infty f_a^{(-k)}(p_0)\right\}
}.
\label{J_a}
\end{equation}

As $a$ becomes small every point of $\overline{{\rm DSPP}}$ moves around. At $a=0$ limit the Julia set must disappear and suddenly appear IVPPs different for each period. To understand what happens it will be useful to study simple examples and see their pictures. 
Let us consider two rational maps of 2 and 3 dimensions:
\begin{equation}
(x,y)\to f_a(x,y)=\left(x{1-y\over 1-x-a},\ y{1-x\over 1-y}\right),
\label{2dMoebius}
\end{equation}
\begin{equation}
(x,y,z)\to f_{a,b}(x,y,z)=\left(x{1-y+yz\over 1+a-z+zx},\ y{1+b-z+zx\over 1-x+xy},\ z{1-x+xy\over 1-y+yz}\right),
\label{deformed 3dLV}
\end{equation}
where we used the notation $(x,y,z):=(x_1,x_2,x_3)$, for convenience.

When $a=0$, these maps are completely integrable. In fact the map (\ref{2dMoebius}) has an invariant $r=xy$ at $a=0$ and the map becomes a M\"obius map, while (\ref{deformed 3dLV}) has two invariants $r=xyz$ and $s=(1-x)(1-y)(1-z)$ at $a=b=0$ and becomes the 3 dimensional Lotka-Volterra map (3dLV), which is known integrable.

Since we are interested in how the periodic points move as $a$ changes, it will be convenient if we can see explicitly paths of periodic points of any particular periods as the parameter  changes. After elimination of the parameter from the periodicity conditions (\ref{ppc}) in the period 2 case of the map (\ref{2dMoebius}), for example, we find an algebraic curve
\begin{equation}
(2-x)^2(1-y)^2-3(1-x)(1-y)(2-x)
+(1-x)^2(2-x+xy)=0,
\label{G2}
\end{equation}
which is drawn by red in Fig.1a. Similar curves we can draw for the period 3(blue) and period 4(green) cases, but the corresponding formulae become more complicated. In summary we have found:
\begin{enumerate}
\item
All points of the Julia set approach $(x,y)=(1,1)$ along the closed curves crossing at $(1,1)$.
\item
All stable periodic points approach certain points on open curves different for each period. 
\end{enumerate}

\begin{minipage}{7cm}
\vspace{-0cm}

\begin{center}
\includegraphics[scale=0.35]{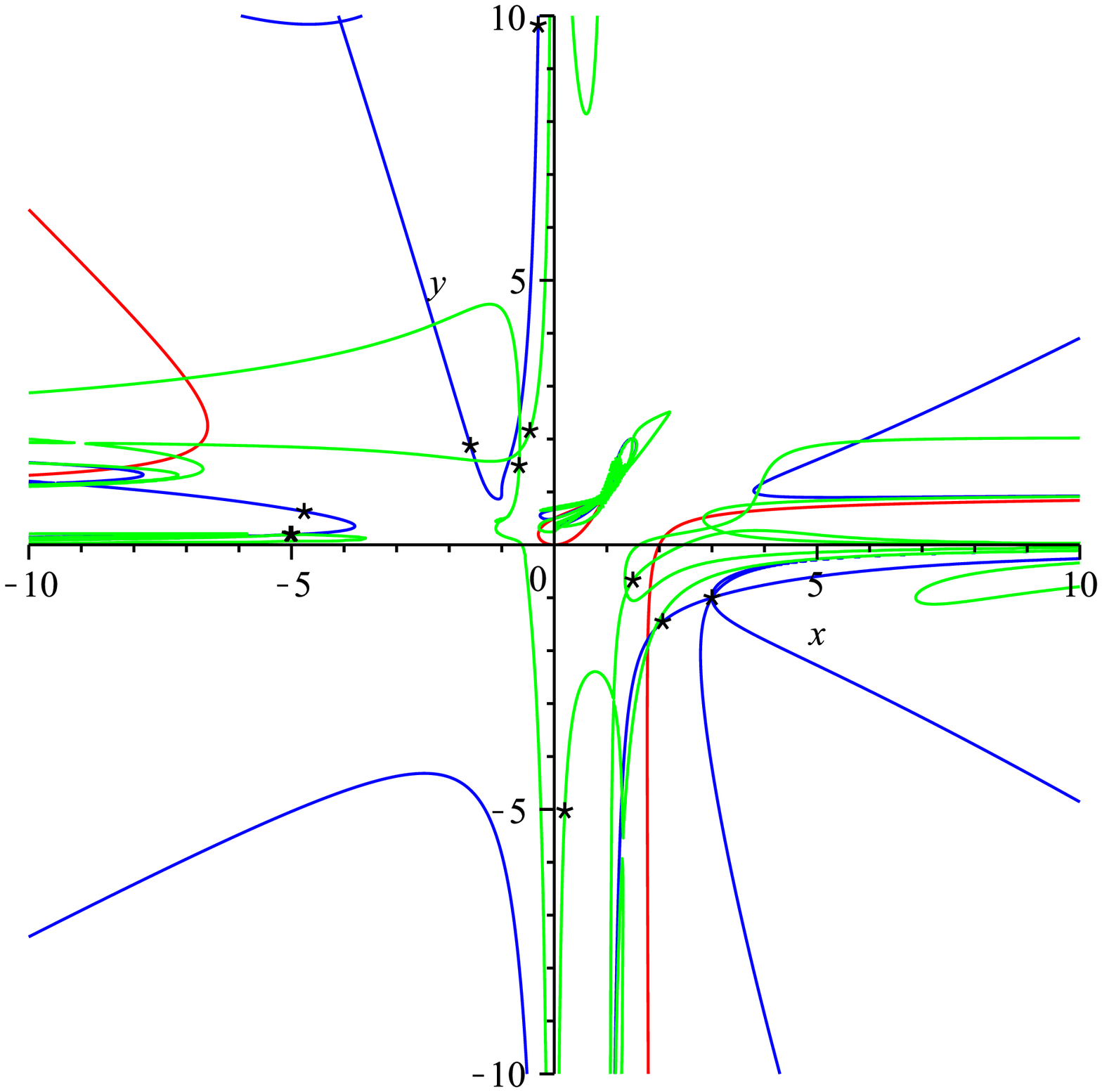}
\vspace{0cm}

Fig.1a\quad period 2,3,4 curves
\end{center}
\end{minipage}
\qquad
\begin{minipage}{7cm}
\vspace{-0cm}

\begin{center}
\includegraphics[scale=0.35]{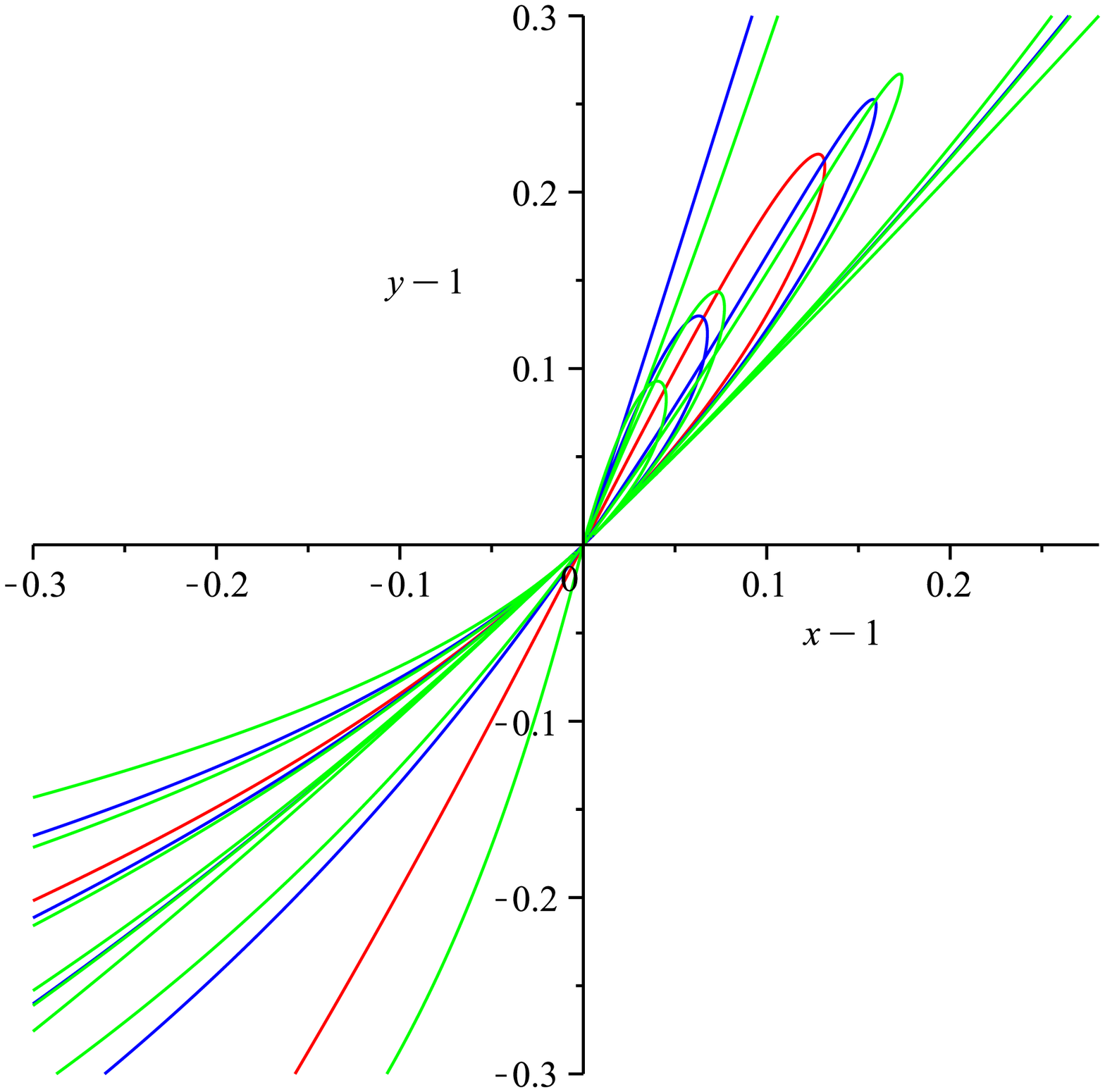}
\vspace{0cm}

Fig.1b\quad details near $(x,y)=(1,1)$
\end{center}
\end{minipage}
\bigskip

The points where the stable points of period 3 and period 4 approach in the $a=0$ limit are shown by asterisks in Fig.1a. The stable period 2 points move to infinity in this limit. \\

We can apply the same method to see how the periodic points of the deformed 3dLV map (\ref{deformed 3dLV}) move. In this case we have two parameters $a,b$. As they change every periodic point must move on a surface instead of a curve. In fact, after elimination of $a,b$ from the period 2 conditions, we obtain a formula
\begin{equation}
K^{(2)}(x,y,z)=0,
\label{K}
\end{equation}
which represents a complex surface in $\mathbb{C}^3$. Since the formula (\ref{K}) is too big, we present it in Appendix. The surface is shown in Fig.2a. It should be compared with the red curve of Fig.1 for the map (\ref{2dMoebius}). The singular locus at $(1,1)$ of Fig.1b corresponds to the curves along which the surface intersect with itself. We show a detail of one of the intersection curves in Fig.2b. 

\begin{minipage}{7cm}
\vspace{-0cm}

\begin{center}
\includegraphics[scale=0.35]{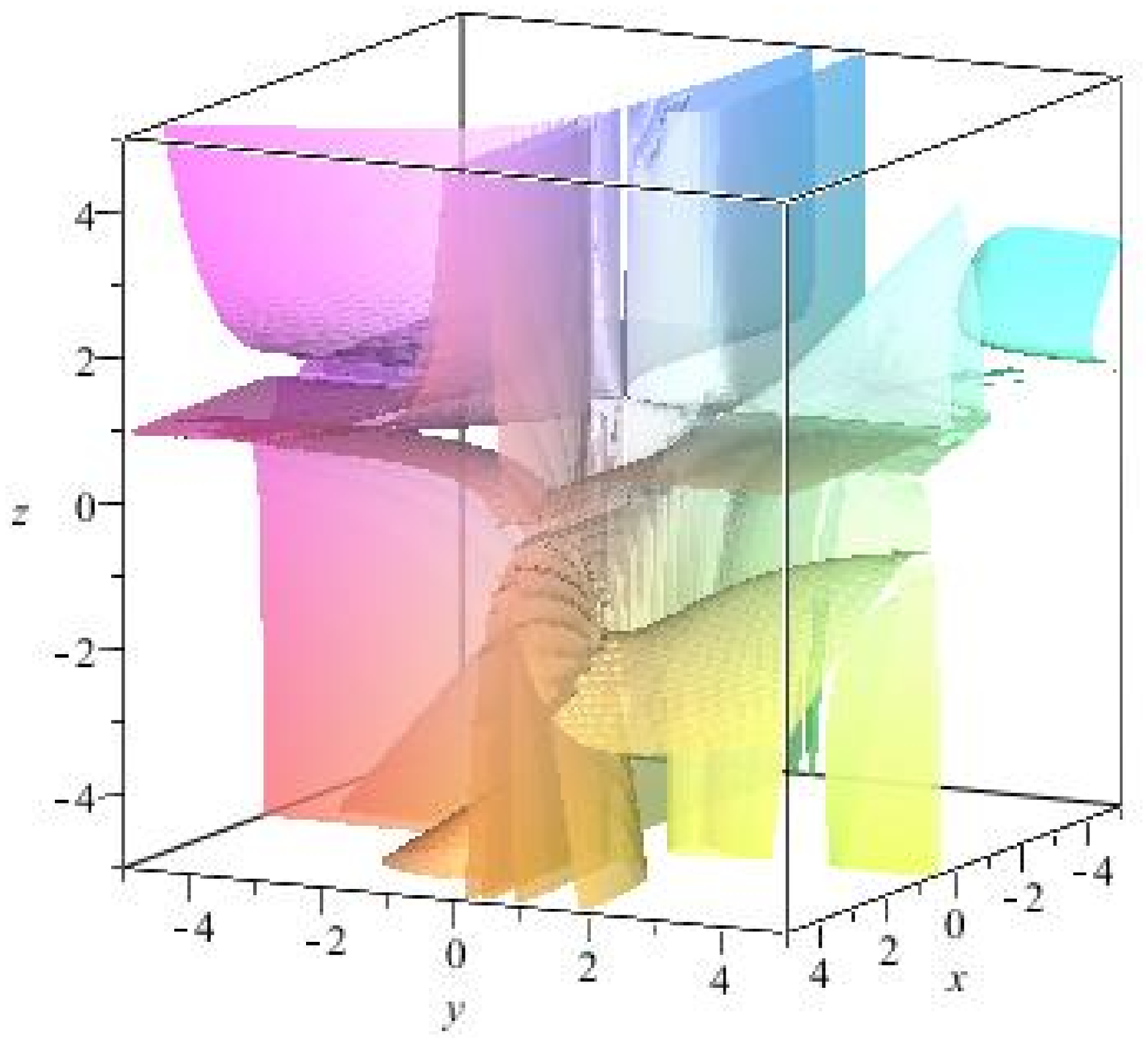}
\vspace{0cm}

Fig.2a\quad period 2 surface $K^{(2)}=0$
\end{center}
\end{minipage}\qquad\quad
\begin{minipage}{7cm}
\vspace{-0cm}

\begin{center}
\includegraphics[scale=0.35]{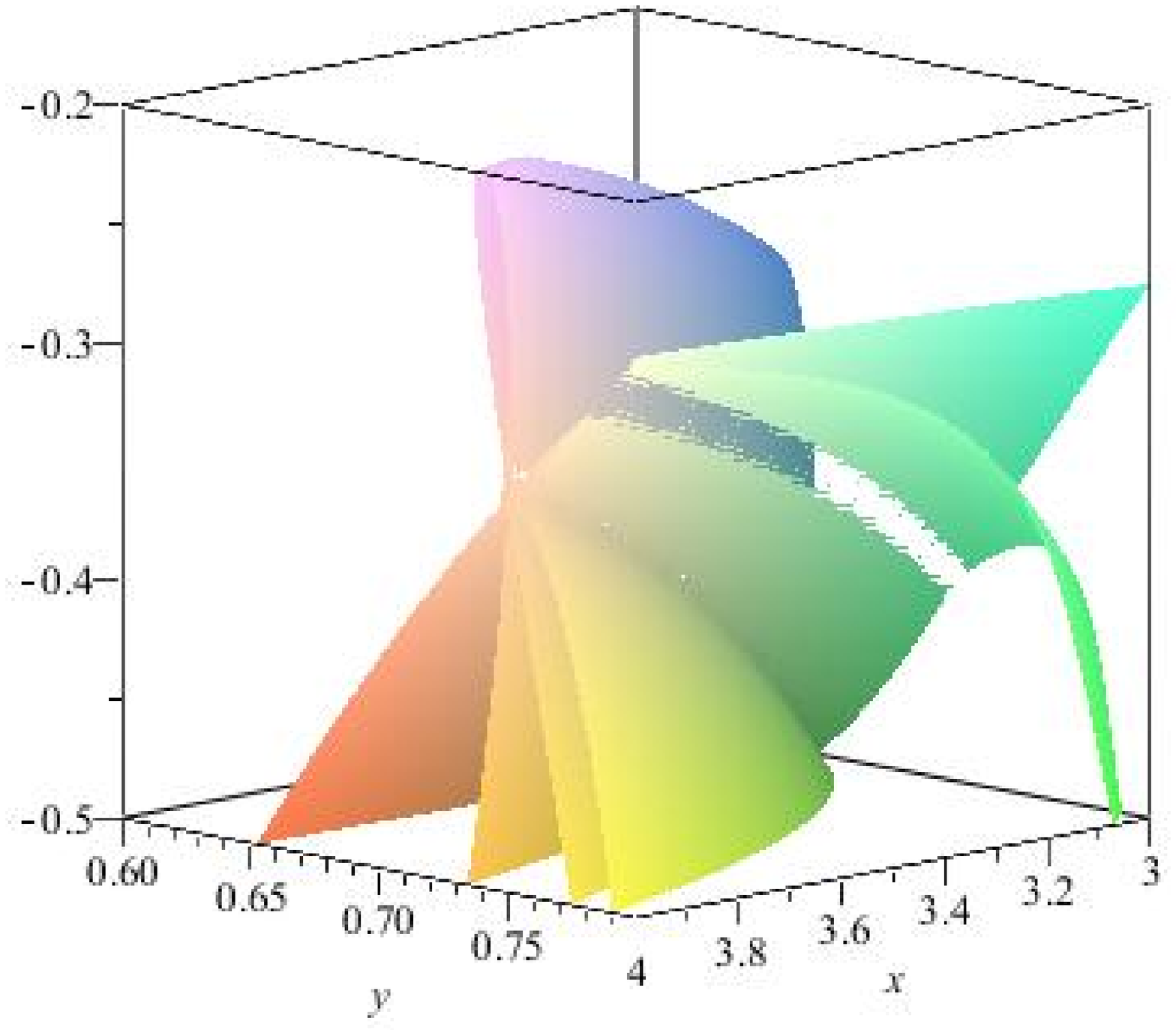}
\vspace{0cm}

Fig.2b\quad details near $\Lambda^+$
\end{center}
\end{minipage}
\bigskip

It is not difficult to find the curves of selfintersection of the surface (\ref{K}). They consist of two independent pieces
\begin{equation}
\Lambda^+:=\left(1-{1\over t}, {1\over 1-t},t\right),\quad 
\Lambda^-:=\left( {1\over 1-t}, 1-{1\over t},t\right),\qquad
t\in\mathbb{C}.
\label{indeterminate curve}
\end{equation}
The curve of Fig.2b corresponds to $\Lambda^+$.  Although the surfaces of large periodic numbers become more complicated we can show that all surfaces of different periods intersect along the same curves $\Lambda^\pm$. 

Once again we find that the points of the Julia set approach $\Lambda^\pm$ as $a,b$ become small. On the other hand, all stable points go somewhere else on the surfaces. According to the IVPP theorem they must go to IVPPs, which we discuss in the following.\\

When the control parameters vanish the maps become integrable and all periodic points are on IVPPs. We can derive IVPPs of all periods iteratively using the method developed in \cite{SS5, YSW, SYW}. For example we find $\gamma$'s of  (\ref{gamma}) as
\begin{eqnarray}
\gamma^{(3)}&=&r+3,\nonumber\\
\gamma^{(4)}&=&r+1,\nonumber\\
\gamma^{(5)}&=&r^2+10r+5,\nonumber\\
 {\rm etc.}&&
\label{2dMobiusgamma}
\end{eqnarray}
from the map (\ref{2dMoebius})  with $a=0$ and  
\begin{eqnarray*}
\gamma^{(2)}&=&s+1,\\
\gamma^{(3)} &=&(s-r)^2+(r+1)(s+1),\\
\gamma^{(4)}&=&(s-r)^3+s(r+1)^3,\\
{\rm etc.},&&
\end{eqnarray*}
from the map (\ref{deformed 3dLV}) with $a=b=0$. 

From information (\ref{2dMobiusgamma}) we can write down the IVPPs of all periods as follows:
\begin{equation}
v^{(n)}=\Big\{ x,y\in \mathbb{C}^2 \Big|\ xy+\tan^2(\pi m/n)=0,\ 
m=1,2,..,n-1\Big\},\qquad  n=3,4,5,...
\label{2dMoebiusIVPPs}
\end{equation}
in the case of the map (\ref{2dMoebius})  with $a=0$. They are hyperbolic curves as drawn in Fig.3.
\vspace{-0cm}

\begin{center}
\includegraphics[scale=0.3]{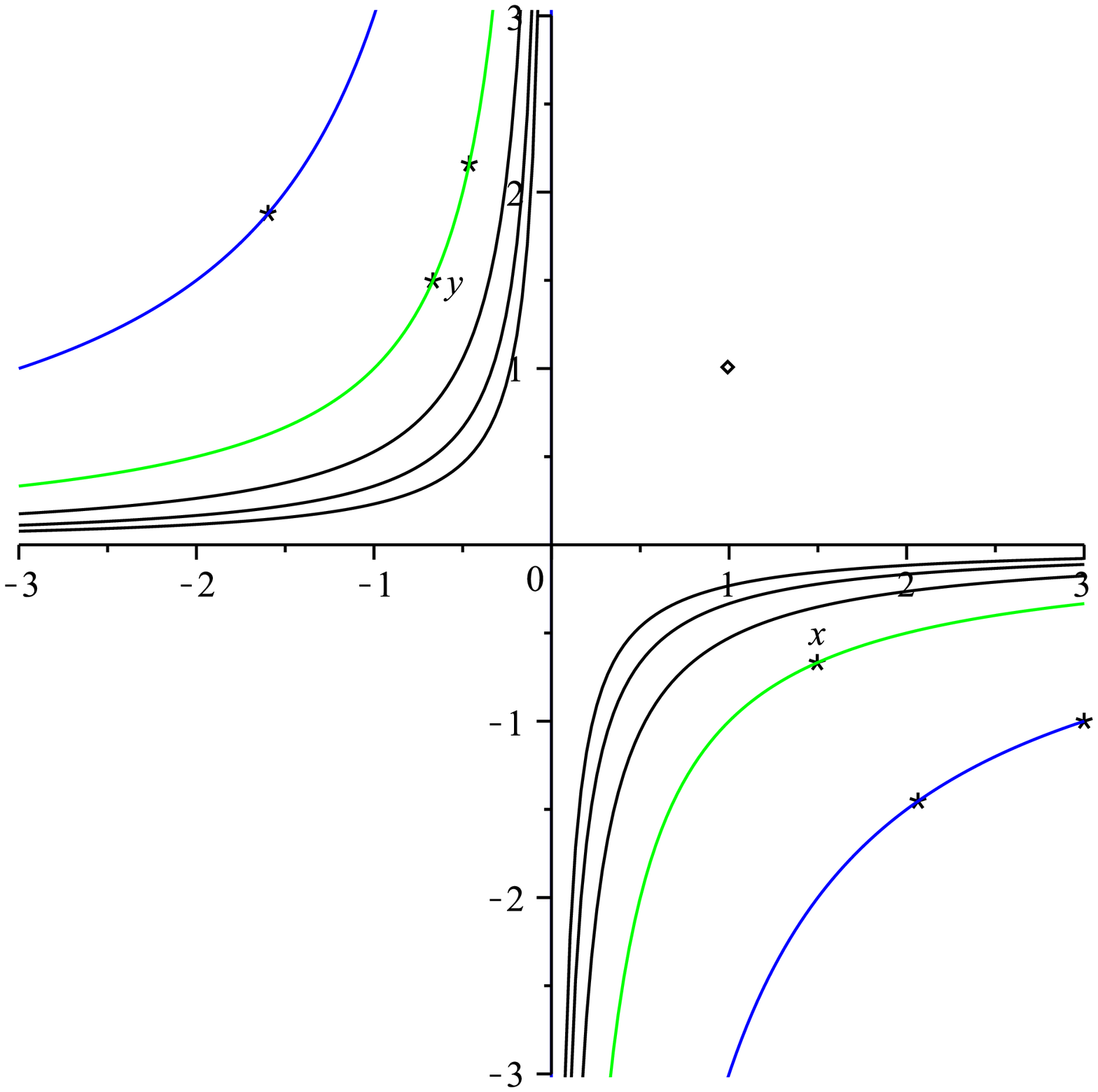}
\vspace{0cm}

Fig.3\quad  IVPPs of (\ref{2dMoebiusIVPPs})
\end{center}
\bigskip

In the 3dLV case the IVPPs in the $(x,y,z)$ space become much more complicated to show in a picture. We present only those of period 2 and period 4 in Fig.4a. Here the red surface is the IVPP of period 2, while the green one is of period 4. For higher period cases it is more convenient to show them in the $(r,s)$ space. In Fig.4b we show the curves of $\gamma^{(n)}=0$ in the $(r,s)$ space for $n=2$ to 10.

\begin{center}
\begin{minipage}{7cm}
\vspace{0cm}

\begin{center}
\includegraphics[scale=0.35]{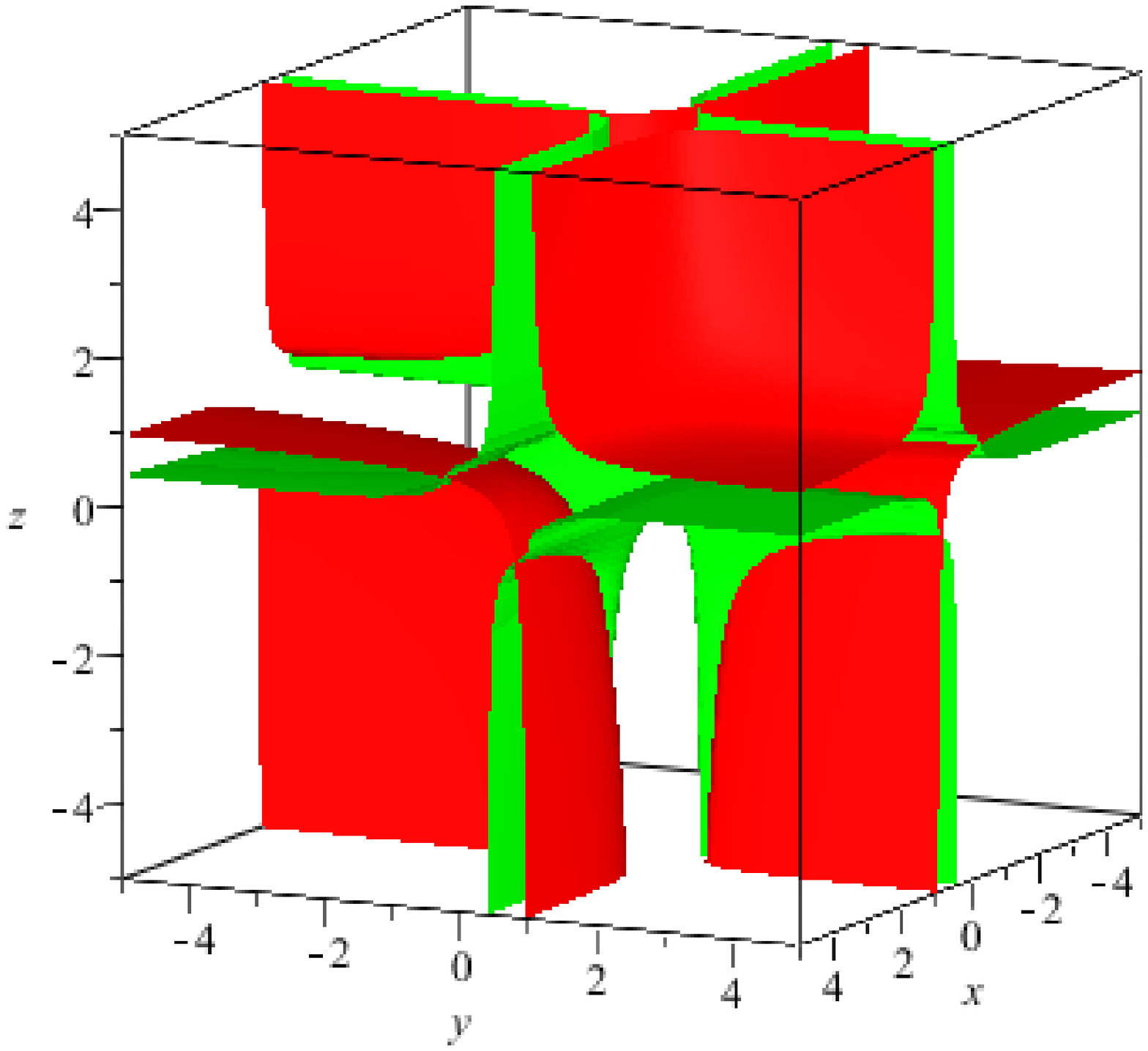}
\vspace{0.5cm}

Fig.4a\quad  period 2 and 4 IVPPs of 3dLV
\end{center}
\end{minipage}
\qquad
\begin{minipage}{8cm}
\vspace{1cm}
\begin{center}
\includegraphics[scale=0.4, bb=0 0 564 400]{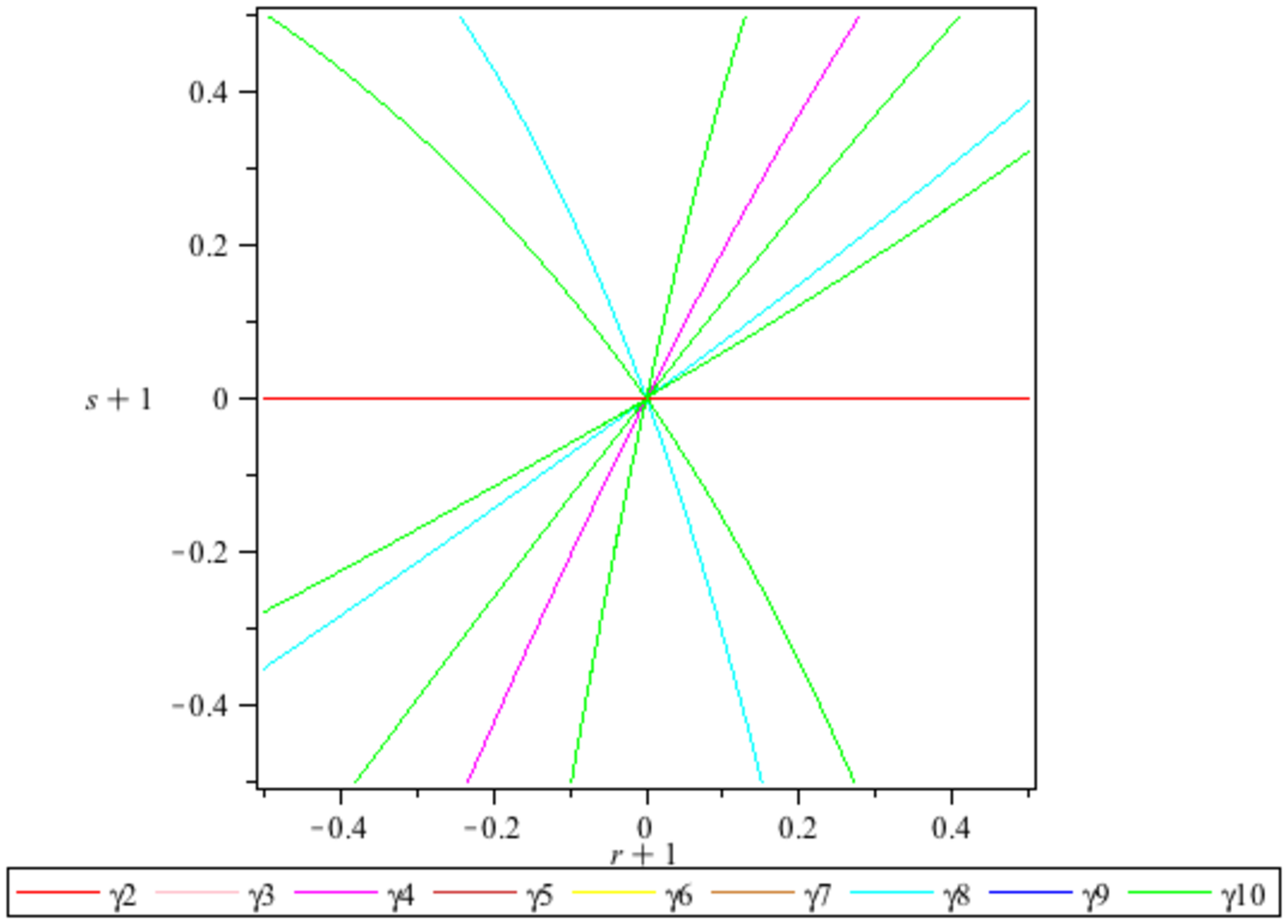}
\vspace{0.2cm}

Fig.4b\quad IVPPs of 3dLV in $(r,s)$
\end{center}
\end{minipage}
\end{center}
\bigskip

We notice that all IVPPs intersect at $(r,s)=(-1,-1)$. The curves for $n=3,6,9$ do not appear in this figure because they are not on the real $(r,s)$ plane, beside the intersection point. Moreover we can convince ourselves that this point of intersection corresponds exactly the curves $\Lambda^\pm$, which we discussed above. They are nothing but the varieties formed by indeterminate points of the map. Since every point on $\Lambda^\pm$ are occupied by periodic points of all periods we called them varieties of singular points (VSP) in our previous papers\cite{SS1,IVPP,SS5}.

%%%%%%%%%%%%%%%%%%%%%%%%%%%%%%%%%%%%%%%%%%%%%%%
\section{The fixed point theorem at the transition point}

We reviewed in \S2 the phenomenon which takes place when a dynamical map transits from a nonintegrable phase to an integrable phase. The purpose of this section is to associate this phenomenon with the transition of the fixed point theorems, which we discussed in \S1. Throughout this section we again consider iteration of a rational map of $d$ dimension, which has the control parameter $a$.

\subsection{Nonintegrable phase}

When $a\ne 0$ the periodic points of period $n$ form a subset of DSPP consisting of a finite number of discrete points. If the dimension of this set is $M_n$ we can label the points by an integer $m\le M_n$, such that every element of DSPP is labeled by two integers $(m,n)\in \mathbb{N}\times\mathbb{N}$.

The position of the point $(m,n)$ in $\mathbb{C}^d$ can be specified by the function 
\[
x: (m,n)\to x(m,n)\in {\rm DSPP}.
\]
We identify this function $x(m,n)$ with the function $\phi$ of the Yanofsky diagram (\ref{Yanofsky}), hence the set $Y$ is the DSPP in this case. 

For the purpose of considering the diagonal argument we define a function $\alpha$ such that
\[
\alpha: x(m,n)\to x(m+1,n), \quad \mod(M_n).
\]
Apparently $\alpha$(DSPP)=DSPP. Since $M_n\ge 2$ for all nonlinear maps, this $\alpha$ does not have a fixed point. 
%\begin{equation}\begin{array}{ccc}\mathbb{N}\times \mathbb{N}&\maprightu{\phi}&{\rm DSPP}\cr\mapupl{\Delta}&&\mapdownr{\alpha}\cr\mathbb{N}&\maprightd{\psi}&{\rm DSPP}\cr\end{array}\label{YanofskyLT1}\end{equation}
The existence of $\alpha$ which has no fixed point tells us, owing to the theorem LT1, that there is a function $\psi(-)$ which is not representable by a function $\phi(-,n)$ of any $n$.\\

Now let us consider the meaning of this result. We associated nonintegrable maps with LT1 of  Lawvere's theorem, which can explain various paradoxes in logic and mathematics. Indeed, despite of our assumption of a map being deterministic, the absence of the fixed point yields the Julia set, which makes some of orbits unpredictable. 

It will be useful to compare our argument with Cantor's diagonal argument on the cardinal numbers, which shows that the countable numbers are not sufficient to count real numbers. In our problem we consider periodic points consisting of countable number of subsets. But the fractal nature of the Julia set makes their closure dense, such that its cardinal number is one of continuum. Since the Julia set is a subset of {$\overline{{\rm DSPP}}$}, the cardinal number of {$\overline{{\rm DSPP}}$} itself  becomes one of continuum. 
 Schematically we have the following correspondence:
\begin{equation}
\begin{array}{ccccccc}
{\rm unpredictable}&\Leftarrow&\exists\ {\overline{{\rm DSPP}}}&\Leftarrow 
\nexists\ {\rm fixed\ point}&\Leftarrow&{{\rm DSPP}}\\
{\rm uncountable}&\Leftarrow&\exists\ {\rm continuum}&\Leftarrow 
\nexists\ {\rm fixed\ point}&\Leftarrow&{\rm countable}\\
\end{array}
\label{unpredictable}
\end{equation}

It is also interesting to notice that, in the symbolic dynamics, every chaotic orbit can be associated with an irrational number\cite{Devaney}.

\subsection{Integrable {phase}}

We now consider the same map but at $a=0$. In this case, as we discussed in \S1, we are to consider Brouwer's fixed point theorem, since the objects which we concern are IVPPs. 

Since the map in this theorem is continuous, it will be convenient to introduce a real parameter $t\in [t_0,t_1]$, and consider explicitly a continuous map ${\rho(t_1,t_0):} D(t_0)\to D(t_1)$ as $t$ changes from $t_0$ to $t_1$. We can think of this map as a propagation of the disk $D(t_0)$ to $D(t_1)$, or equivalently
\[
{Cyl(\rho,D):=\Big\{ D(t):=\rho(t,t_0)\big(D(t_0)\big) \sim D(t_0)\Big| 
\rho(t',t): D(t)\to D(t'),\ t,t'\in [t_0,t_1]\Big\},}
\]
%$\{D(t)|t_0\le t\le t_1\}$, 
which we call a `cylinder'. Then Brouwer's fixed point theorem is equivalent to the statement: 
\bigskip

\noindent
{\bf Corollary}

{\it Any two cross sections of the `cylinder' have a fixed point 
{with respect to the map $\rho$}.}
\bigskip

{We want to apply this theorem to our case. Since our objects are IVPPs,} we generalize the disk to topological space $X$ and consider a `cylinder' 
\begin{equation}
{Cyl(\rho,X):=\Big\{ X(t):=\rho(t,t_0)\big(X(t_0)\big) \sim X(t_0)\Big| 
\rho(t',t): X(t)\to X(t'),\ t,t'\in [t_0,t_1]\Big\},}
\label{Cyl(rho,X)}
\end{equation}
%\[
%Cyl:=\{X(t)|\ t_0\le t\le t_1\},
%\]
which is obtained by a continuous map ${\rho(t,t_0):} X(t_0)\to X(t)$ on the base space $X_0:=X(t_0)$. This `cylinder' defines a fiber bundle $\pi: Cyl\to X_0$. In particular we call $X(t)$ a cross section of the `cylinder' at $t$. A diagonal cross section $D_c$ is defined by the cross section which includes at least a point of $X(t_0)$ and a point of $X(t_1)$,
 {\it i.e.,} $D_c\cap X(t_0)\ne \emptyset,\ D_c\cap X(t_1)\ne \emptyset$.

{Let us denote by $I(t, D_c)$ the intersection of $X(t)$ with $D_c$.
We say that $Cyl(X, \rho)$ is representable by $D_c$ if $\rho(t, t')$ for all $t, t'$ is obtained by extending an injection of a subset of $I(t, D_c)$.}
%We say that a diagonal cross section $D_c$ of the `cylinder' is representable if there is a $t\in [t_0,t_1]$ such that $D_c$ is a continuous map of the cross section $X(t)$.
The next statement follows immediately to the Corollary.
\bigskip

\noindent
{\bf Proposition}

{\it If an arbitrary diagonal cross section {$D_c$} of the `cylinder' $Cyl$ is representable, any two cross sections have a fixed point.}\\

{We are now ready to apply this result to our integrable maps. We want to identify $X(t)$ with an IVPP. But it is not obvious because $t$ of $X(t)$ is continuous while an IVPP is defined only for every integer. We must either restrict $t$ to integers or extend an IVPP $v^{(n)}$ to $v^{(t)}$. We adopt the latter, because, as far as we consider integrable maps, we know explicitly the map $f^{(n)}(x)$ as well as $v^{(n)}$ as functions of $n$ in general.} For example from the expression of $v^{(n)}$ of  (\ref{2dMoebiusIVPPs}) we immediately find
\begin{equation}
v^{(t)}=\{ x,y\in \mathbb{C}^2 |\ xy+\tan^2(\pi /t)=0\},\qquad t\in\mathbb{R}
\label{continuous IVPP}
\end{equation}
as the continuous version of IVPPs. {When $t$ is irrational $v^{(t)}$ fixes an invariant variety filled by quasi periodic orbits alone specified by $t$. This result agrees with other cases we studied in} \cite{SS4}.

In our previous works\cite{SS1, IVPP, SS2, SS3} many integrable discrete equations, such as the Hirota-Miwa equation\cite{H, M}, the Quispel-Roberts-Thompson equation\cite{QRT}, discrete Euler top, etc., are shown that  periodic points are on IVPPs. All possible continuous limits of these equations are also integrable. From the single Hirota-Miwa equation, for example, infinitely many soliton equations which belong the KP hierarchy can be derived. Their solutions are given by a universal function, called the $\tau$ function, irrespective to discrete or continuous.

Based on these observations we want to conjecture the following:

\noindent
{\bf Conjecture}

{\it If $f^{(n)}(x), \ n\in\mathbb{Z}$ is a sequence of discrete integrable maps for all $x$, there exists a continuous integrable map $f^{(t)}(x),\ t\in \mathbb{R}$.}
\bigskip

This enables us to identify the continuous map $\rho(t',t)$ of the cylinder (\ref{Cyl(rho,X)}) with the map of the projective resolution of the triangulated category and $v^{(n)}$ with the cross section $X(t)$ at $t=n$.  The `cylinder' $Cyl$ must be extended to $t\in [2,\infty)$. Then, according to Proposition, $v^{(n)}$ must intersect with $v^{(n+1)}$ at some point on $v^{(n)}$. It is a fixed point of all IVPPs.\\

Let us summarize our result of this subsection by means of the Yanofsky diagram\cite{Yanofsky}. 
In the diagram the function $\phi$ maps the cylinder $\mathbb{R}\times X$ to the `cylinder' $Cyl\in \mathbb{C}^d$, while the function $\psi$ is the continuous map which maps $X_0$ to the diagonal cross section $D_c\in Cyl$. At the same time every cross section of the `cylinder' can be mapped by the function $\alpha$ continuously to any other cross sections including the diagonal cross section $D_c$. This explanation is certainly not sufficient to prove our proposition from mathematical point of view.  Nevertheless we would like to emphasize that this is always true as far as we studied many nonlinear maps with IVPPs in the integrable regime.  
 
Before closing this section it will be interesting to compare our argument with the proof of G\"odel's first incompleteness theorem.
Integrability of a map means that all orbits are predictable  for arbitrary initial points. The existence of IVPPs guarantees the predictability of the map. But, according to LT2, they also allows for the VSP to exist, where every point becomes indeterminate since periodic points of all periods degenerate there altogether. 

In the proof of G\"odel's first incompleteness theorem one assumes the existence of a formal axiomatic system (FAS) which can decide all possible propositions from a finite number of axioms. But the FAS itself possesses a fixed point which shows the existence of a proposition undecidable whether it is true or false. Namely 
the FAS is incomplete. Thus we have the correspondence, schematically, as follows:
\begin{equation}
\begin{array}{ccccccc}
{\rm predictable}&\Rightarrow&\exists\ {\rm IVPP}&\Rightarrow&
\exists\ {\rm fixed\ point}&\Rightarrow&{\rm indeterminate},\\
{\rm decidable}&\Rightarrow&\exists\ {\rm FAS}&\Rightarrow&
\exists\ {\rm fixed\ point}&\Rightarrow&{\rm incomplete}.\\
\end{array}
\label{incomplete}
\end{equation}

%%%%%%%%%%%%%%%%%%%%%%%%%%%%%%%%%%
\section{Conclusion}

We can summarize our results of \S2 and \S3 as
\begin{equation}
\begin{array}{ccccccccc}
{\rm `unpredictable'}&\Leftarrow&{\rm nonintegrable}&\Leftarrow&\exists\ \overline{{\rm DSPP}}&\Leftarrow&\nexists\ {\rm fixed\ point},&&\\
&&\updownarrow &&
\updownarrow ({\rm I})&&
\updownarrow ({\rm L})&&\\
&&{\rm integrable}&\Rightarrow&\exists\ {\rm IVPP}&\Rightarrow&\exists\ {\rm fixed\ point}&\Rightarrow&{\rm incomplete}.
\end{array}
\label{results}
\end{equation}
Here `unpredictable' means not totally predictable. The transition (I) owes to the IVPP theorem. On the other hand the transition (L) is not clear from our arguments in \S2 and \S3, because the existence of fixed points in integrable regime has been proved by means of Brouwer's fixed point theorem. If the fixed points in (\ref{results}) are those of Lawvere's fixed point theorem, we can say the transition (L) is implied by the contrapositive of Lawvere's fixed point theorem. Therefore we now focus our attension to this problem. Unfortunately we do not know how to extend the conditions for the proposition LT2 to include Brouwer's fixed point theorem. Nevertheless we are able to show their correspondence indirectly as follows. 

First we notice that Brouwer's fixed point theorem was proved\cite{Tanaka} being equivalent to Arrow's impossibility theorem well known in social choice theory\cite{Arrow}. Moreover there exist some arguments\cite{Abramsky} to describe Arrow's impossibility theorem by means of the category theory. Therefore it is quite natural to associate Arrow's impossibility theorem with Lawvere's fixed point theorem.
Mathematical details of this proof will be presented in our forthcoming paper. Based on this argument we are convinced ourselves that the fixed points which appear in integrable maps are those required the existence by means of Lawvere's fixed point theorem.\\

In the rest of this section we would like to discuss the meaning of our results. 
We could associate the transition scheme of integrable nonintegrable maps with the transition (\ref{LT scheme}) of Lawvere's fixed point theorem to its own contrapositive assertion. The fixed point is the VSP in our case. As a nonintegrable map becomes integrable, there appears VSP together with IVPPs, while the Julia set totally disappears. We can look this phenomenon from the opposite side. When an integrable map is perturbed, the fixed point disappears together with IVPPs and the Julia set emerges from the fixed point. Therefore the fixed point is the source of both IVPPs and the Julia set.\\

The self-referential paradoxes, which have been described by Lawvere's fixed point theorem, are those in logic, mathematics or in computer programming. From the correspondence (\ref{unpredictable}) the unpredictability of the nonintegrable phase is associated with the
appearance of nonrepresentable objects. As the phase changes to the integrable one the system becomes, according to (\ref{incomplete}), incomplete. Therefore we must conclude that a system of dynamical map is either incomplete or not totally predictable in the sense of (\ref{results}).

Since the problems we have discussed in this paper are concerned with phenomena in physics, there arises a question, whether the unpredictability and/or incompleteness is the nature of the phenomena  or their description?
When one obtains two different answers which contradict with each other, the paradox is usually thought being due to incompleteness of the logic and/or language. We will adopt this view point because the description of nonlinear phenomena we have discussed is a theory of classical physics. We know that the classical theory of physics is not sufficient to explain real physical phenomena. The best we can conclude from our discussion in this paper is that the classical description of nonlinear dynamics is neither complete nor totally predictable.

\appendix

\section{Period 2 surface (\ref{K}) after elimination of $(a,b)$}

\begin{eqnarray*}
&&K^{(2)}(x,y,z)\\
&&=3y^5-2y^2z^4x^3-134z^3x^4y^3+24z^3x^4y^2+358z^2x^4y^3
-3z^3x^4y-102z^2x^4y^2+9z^2x^4y
\\
&&-3y^3z^5x^2+13y^4z^5x+y^3z^5x-2y^4z^5x^2-7y^4-5y^6
+y^7-129x^4y^4+117x^3y^4-x^7y^7
\\&&+146x^4y^3-101x^4y^2+40x^4y-6y^4z^5
+28xy^4+4zx^4-83y^4x^2-y^7x-y^7z^5+5y^7z^4
\\&&
-10y^7z^3+10y^7z^2-5y^7z-477zx^3y^4
-591x^4y^4z^2-385zx^4y^3+183zx^4y^2-45zx^4y
\\&&+462x^4y^4z+x^7y^7z^3-5x^7y^6z^3+10x^7y^5z^3
-10x^7y^4z^3+5x^7y^3z^3-x^7y^2z^3+2x^6y^7z^4
\\&&-5x^6y^7z^3+5x^6y^7z^2-8x^6y^6z^4+26x^6y^6z^3-37x^6y^6z^2
+13x^6y^5z^4-55x^6y^5z^3-x^6y^2z^4
\\&&-11x^6y^4z^4+60x^6y^4z^3-154x^6y^4z^2+5x^6y^3z^4-35x^6y^3z^3
+121x^6y^3z^2+106x^6y^5z^2
\\&&+10x^6y^2z^3-49x^6y^2z^2-x^6yz^3+8x^6yz^2+x^5y^8z^4-3x^5y^8z^3+3x^5y^8z^2
+x^5y^7z^5
\\&&-10x^5y^7z^4+30x^5y^7z^3-35x^5y^7z^2
-2x^5y^6z^5+31x^5y^6z^4-107x^5y^6z^3+140x^5y^6z^2
\\&&+x^5y^5z^5-42x^5y^5z^4+186x^5y^5z^3-289x^5y^5z^2
+28x^5y^4z^4-174x^5y^4z^3+345x^5y^4z^2
\\&&-10x^5y^3z^4+89x^5y^3z^3-239x^5y^3z^2+2x^5y^2z^4
-24x^5y^2z^3+88x^5y^2z^2+3x^5yz^3+9x^7yz
\\&&-13x^5yz^2+x^4y^8z^5-2x^4y^8z^4+x^4y^8z^3
-9x^4y^7z^5+30x^4y^7z^4-43x^4y^7z^3-42x^3y^7z^2
\\&&+35x^4y^7z^2-15x^4y^7z+x^4y^6z^6-215x^4y^6z^2
+21x^4y^6z^5-97x^4y^6z^4+198x^4y^6z^3-x^7z
\\&&+111x^4y^6z-2x^4y^5z^6-18x^4y^5z^5+127x^4y^5z^4
-366x^4y^5z^3+506x^4y^5z^2-315x^4y^5z
\\
&&+x^4y^4z^6+6x^4y^4z^5-x^4y^3z^5-2x^3y^8z^5+6x^3y^8z^4-6x^3y^8z^3+2x^3y^8z^2
+18x^3y^7z^5
\\&&-60x^3y^7z^4+74x^3y^7z^3+12x^3y^7z-4x^3y^6z^6
-40x^3y^6z^5+188x^3y^6z^4-298x^3y^6z^3
\\&&+230x^3y^6z^2-92x^3y^6z+8x^3y^5z^6-75x^7y^4z
+33x^3y^5z^5+296x^3y^5z-4x^3y^4z^6+y^4z^6
\\&&+54x^2y^7z^4-80x^2y^7z^3+56x^2y^7z^2-18x^2y^7z
+2x^2y^7+6x^2y^6z^6+23x^2y^6z^5-159x^2y^6z^4
\\&&+288x^2y^6z^3-236x^2y^6z^2+93x^2y^6z-15x^2y^6
-12x^2y^5z^6-7x^2y^5z^5+47x^2y^5+6x^2y^4z^6
\\&&+5xy^7z^5-21xy^7z^4+34xy^7z^3-26xy^7z^2+9xy^7z
-4xy^6z^6+52xy^6z^4-122xy^6z^3+9y^5
\\
&&+114xy^6z^2-46xy^6z+6xy^6+8xy^5z^6-18xy^5z^5
-41xy^5z^4+172xy^5z^3-206xy^5z^2
\\&&+102xy^5z-17xy^5-4xy^4z^6+y^6z^6
-2y^6z^5-7y^6z^4+28y^6z^3-37y^6z^2+22y^6z-2y^5z^6
\\&&+9y^5z^5-9y^5z^4-16y^5z^3+42y^5z^2-33y^5z
+y^8x^2z^5-4y^8x^2z^4+6y^8x^2z^3-14x^2y^7z^5
\\
&&-4y^8x^2z^2-50x^5y+2x^4y^7
+y^8x^2z-7x^4+9x^5-5x^6-59y^5x^3+68x^4y^5+115x^5y^4
\\&&-150x^5y^3+117x^5y^2-19x^4y^6
-2x^3y^7+16x^3y^6-54x^5y^5-11x^6y^6+45x^6y^5-33x^7y^2z
\\&&-95x^6y^4+115x^6y^3-81x^6y^2+31x^6y+15x^5y^6
+7x^7y^6-21x^7y^5+35x^7y^4-35x^7y^3
\\&&+21x^7y^2-7x^7y+x^6y^7-2x^5y^7+198x^5y^5z
+308zx^5y^3-179zx^5y^2+54zx^5y-40x^7y^5z^2
\\&&-313x^5y^4z-6zx^5-77zx^5y^6
+16x^5y^7z-3x^7y^7z^2+17x^7y^6z^2+13x^7y^2z^2-2x^7yz^2
\\&&+50x^7y^4z^2-35x^7y^3z^2-x^5y^8z-35x^6yz+3x^7y^7z-19x^7y^6z+51x^7y^5z+65x^7y^3z.
\end{eqnarray*}

\newpage

%%%%%%%%%%%%%%%%%%%%%%%%%%%%%%%%%%%%%%%%%

\end{document}